%
%
%
%
%
%
%
%
%
\newcount\VOL\VOL=3\newcount\YEAR\YEAR=1998
\def\firstpage{1}\def\lastpage{1000}
\def\received{}\def\revised{}
\def\communicated{}
\nopagenumbers
\hsize=13truecm\hoffset=1.43truecm
\vsize=20truecm\voffset=1.9truecm
\font\eightbf=cmbx8
\font\eightit=cmti8
\font\eightrm=cmr8
\font\caps=cmcsc10                    
\font\Caps=cccsc10 scaled \magstep1   

%
\pageno=\firstpage\def\folio{\rm\number\pageno}
\def\DocMath{}
\def\makeheadline{
    \vbox to 0pt{\vskip-40pt\line{\vbox to8.5pt{}\the\headline}\vss}
    \nointerlineskip}
\headline={\ifnum\pageno=\firstpage{\DocMath\hfill\llap{\folio}}
           \else{\ifodd\pageno\rightheadline\else\leftheadline\fi}\fi}
\def\rightheadline{\caps \hfill \rightheadtext \hfill \llap{\folio}}
\def\leftheadline {\caps \rlap{\folio} \hfill \leftheadtext \hfill }
\def\leftheadtext{\ifnum\pageno>\lastpage\else\SAuthor\fi}
\def\rightheadtext{\STitle}

\def\TSkip{\bigskip}
\newbox\TheTitle{\obeylines\gdef\GetTitle #1
\ShortTitle  #2
\SubTitle    #3
\Author      #4
\ShortAuthor #5
\EndTitle
{\setbox\TheTitle=\vbox{\baselineskip=20pt\let\par=\cr\obeylines%
\halign{\centerline{\Caps##}\cr\noalign{\medskip}\cr#1\cr}}%
	\copy\TheTitle\TSkip\TSkip%
\def\next{#2}\ifx\next\empty\gdef\STitle{#1}\else\gdef\STitle{#2}\fi%
\def\next{#3}\ifx\next\empty%
    \else\setbox\TheTitle=\vbox{\baselineskip=20pt\let\par=\cr\obeylines%
    \halign{\centerline{\caps##} #3\cr}}\copy\TheTitle\TSkip\TSkip\fi%
\centerline{\caps #4}\TSkip\TSkip%
\def\next{#5}\ifx\next\empty\gdef\SAuthor{#4}\else\gdef\SAuthor{#5}\fi%
\ifx\received\empty\relax
    \else\centerline{\eightrm Received: \received}\fi%
\ifx\revised\empty\TSkip%
    \else\centerline{\eightrm Revised: \revised}\TSkip\fi%
\ifx\communicated\empty\relax
    \else\centerline{\eightrm Communicated by \communicated}\fi\TSkip\TSkip%
\catcode'015=5}}\def\Title{\obeylines\GetTitle}
\def\Abstract{\begingroup\narrower
    \parskip=\medskipamount\parindent=0pt{\caps Abstract. }}
\def\EndAbstract{\par\endgroup\TSkip\TSkip}
\newbox\TheAdd\def\Addresses{\vfill\copy\TheAdd\vfill
    \ifodd\number\lastpage\vfill\eject\phantom{.}\vfill\eject\fi}
{\obeylines\gdef\GetAddress #1
\Address #2 
\Address #3
\Address #4
\EndAddress
{\def\xs{6truecm}
\setbox0=\vtop{{\obeylines\hsize=\xs#1}}\def\next{#2}
\ifx\next\empty 
     \setbox\TheAdd=\hbox to\hsize{\hfill\copy0\hfill}
\else\setbox1=\vtop{{\obeylines\hsize=\xs#2}}\def\next{#3}
\ifx\next\empty 
     \setbox\TheAdd=\hbox to\hsize{\hfill\copy0\hfill\copy1\hfill}
\else\setbox2=\vtop{{\obeylines\hsize=\xs#3}}\def\next{#4}
\ifx\next\empty\ 
     \setbox\TheAdd=\vtop{\hbox to\hsize{\hfill\copy0\hfill\copy1\hfill}
	        \vskip20pt\hbox to\hsize{\hfill\copy2\hfill}}
\else\setbox3=\vtop{{\obeylines\hsize=\xs#4}}
     \setbox\TheAdd=\vtop{\hbox to\hsize{\hfill\copy0\hfill\copy1\hfill}
	        \vskip20pt\hbox to\hsize{\hfill\copy2\hfill\copy3\hfill}}
\fi\fi\fi\catcode'015=5}}\gdef\Address{\obeylines\GetAddress}

\hfuzz=0.1pt\tolerance=2000\emergencystretch=20pt\overfullrule=5pt

\Title Geometry and analytic theory of Frobenius manifolds
\ShortTitle Frobenius manifolds
\SubTitle
\Author Boris Dubrovin
\ShortAuthor B.Dubrovin
\EndTitle

\Abstract 

Main mathematical applications of Frobenius manifolds are in the theory
of Gromov - Witten invariants, in
singularity
theory, in differential geometry of the orbit spaces of reflection groups
and of their extensions, in the hamiltonian theory of integrable
hierarchies. The theory of Frobenius manifolds establishes remarkable
relationships between these, sometimes rather distant, mathematical
theories.



1991 MS Classification 32G34, 35Q15, 35Q53, 20F55, 53B50
 
\EndAbstract

\Address SISSA, Via Beirut, 2-4, I-34013 TRIESTE, Italy
\Address
\Address
\Address
\EndAddress
%
%
\def\res{\mathop{\rm res}}
\def\rr{{\rm right}}
\def\ll{{\rm left}}
{\bf WDVV equations of associativity} is the problem of
finding of  a
quasihomogeneous, up to at most quadratic polynomial,  function $F(t)$ of
the variables $t=(t^1, \dots, t^n)$
and of a constant nondegenerate symmetric matrix
$\left(\eta^{\alpha\beta}\right)$ such that the following combinations
of the third derivatives
$
c_{\alpha\beta}^\gamma(t):=\eta^{\gamma\epsilon} 
\partial_\epsilon \partial_\alpha \partial_\beta F(t)
$
for any $t$ are structure constants of an asociative algebra
$
A_t ={\rm span}\, (e_1, \dots, e_n) , ~e_\alpha \cdot e_\beta =
c_{\alpha\beta}^\gamma(t) e_\gamma, ~~\alpha, \, \beta =1, \dots, n
$
with the unity $e=e_1$ (summation w.r.t. repeated indices will be
assumed). These equations were discovered by physicists
E.Witten, R.Dijkgraaf, E.Verlinde and H.Verlinde in the beginning of
'90s.
I invented
Frobenius manifolds as the coordinate-free form of WDVV.
%
%
\smallskip
\noindent{\bf 1. Definition of Frobenius manifold} (FM). 

\noindent{\bf 1.1. Frobenius algebra} (over a field $k$; we mainly
consider the
case $k={\bf C}$) is a pair $(A,<~,~>)$, where $A$ is a commutative
associative $k$-algebra with a unity $e$, $<~,~>$ is a symmetric
nondegenerate
{\it invariant} bilinear form $A\times A\to k$, i.e. $<a\cdot
b,c>=<a,b\cdot c>$ for any $a,\, b\in A$. A {\it gradation of the charge
$d$} on $A$ is a $k$-derivation $Q:A\to A$ such that
$<Q(a),b>+<a,Q(b)>=d<a,b>, ~d\in k$. More generally, graded of the charge
$d\in k$ Frobenius algebra $(A,<~,~>)$ over a graded commutative
associative
$k$-algebra $R$ by definition is endowed with two $k$-derivations
$Q_R:R\to R$ and $Q_A:A\to A$ satisfying the properties
$
Q_A(\alpha a) =Q_R(\alpha) a+ \alpha Q_A(a), ~\alpha\in R, ~a\in A
$
$
<Q_A(a), b> + <a,Q_A(b)>-Q_R<a,b> = d <a,b> , ~a,\,
b\in A.
$

\noindent{\bf 1.2. Frobenius structure} of the charge $d$ on the manifold
$M$ is a
structure of a Frobenius algebra on the tangent spaces $T_tM
=(A_t,<~,~>_t)$ depending (smoothly, analytically etc.) on the point $t\in
M$. It must satisfy the following axioms.

{\bf FM1.} The metric $<~,~>_t$ on $M$ is flat (but not necessarily
positive definite). Denote $\nabla$ the Levi-Civita connection for the
metric. The unity vector field $e$ must be covariantly constant, $\nabla
e=0$.

{\bf FM2.} Let $c$ be the 3-tensor $c(u,v,w):=<u\cdot v, w>$, $u,\, v,\,
w\in T_tM$. The 4-tensor $(\nabla_z c)(u,v,w)$ must be symmetric in
$u,\, v,\,  
w, \, z \in T_tM$. 

{\bf FM3.} A linear vector field $E\in Vect(M)$ must be fixed on $M$,
i.e. $\nabla\nabla E=0$, such that the derivations
$Q_{Func(M)}:=E, ~~Q_{Vect(M)}:={\rm id}+{\rm ad}_E
$
introduce in $Vect(M)$ the structure of graded Frobenius algebra of the
given charge $d$ over the graded ring $Func(M)$ of (smooth, analytic etc.)
functions on $M$. We call $E$ {\it Euler vector field}.

Locally, in the flat coordinates $t^1, \dots, t^n$ for the metric
$<~,~>_t$, a FM with diagonalizable (1,1)-tensor $\nabla
E$ is described by a solution $F(t)$ of WDVV associativity equations,
where $\partial_\alpha\partial_\beta\partial_\gamma
F(t)=<\partial_\alpha\cdot \partial_\beta, \partial_\gamma>$, and vice
versa. We will call $F(t)$ {\it the potential} of the FM
(physicists call it {\it primary free energy}; in the setting of quantum 
cohomology it is called Gromov - Witten potential [KM]).


\noindent{\bf 1.3. Deformed flat connection} $\tilde \nabla$ on $M$ is
defined by
the formula
$
\tilde \nabla_u v:= \nabla_u v + z\, u\cdot v.
$
Here $u$, $v$ are two vector fields on $M$, $z$ is the parameter of the
deformation. (In [Gi1] another normalization is used 
$\tilde\nabla \mapsto \hbar \tilde\nabla$, $\hbar=z^{-1}$.) We extend this
to a meromorphic connection on the direct product $M\times {\bf C}$, $z\in
{\bf C}$, by the
formula
$
\tilde\nabla_{d/dz} v=\partial_z v +E\cdot v -z^{-1} \mu\, v ~~{\rm
with}~ \mu:= 1/2(2-d)\cdot {\bf 1} -\nabla E,
$
other covariant derivatives are trivial.
Here $u$, $v$ are tangent vector fields on $M\times {\bf C}$ having zero
components along ${\bf C}\ni z$.
The curvature of $\tilde\nabla$ is equal
to zero. This can be used as a definition of FM [Du3].
So, there locally exist $n$ independent
functions
$\tilde t_1(t;z), \dots, \tilde t_n(t;z)$, $z\neq 0$, such that
$
\tilde\nabla\, d\tilde t_\alpha(t;z)=0, ~\alpha=1, \dots, n
$.
We call these functions {\it deformed flat coordinates}. 
\smallskip 
\noindent{\bf 2. Examples of FMs} appeared first in  2D
topological field theories [W1, W2, DVV].

\noindent{\bf 2.0. Trivial FM:} $M=A_0$ for a graded Frobenius
algebra $A_0$. The potential is a cubic,
$
F_0(t)={1\over 6} <1, (t)^3>, ~~t\in A_0.
$
Nontrivial examples of FM are

\noindent{\bf 2.1. FM with good analytic properties.} They are
analytic perturbations of the
cubic. That means that, in an appropriate system of flat coordinates
$t=(t', t'')$, where all the components of $t'$ have $Lie_Et'={\rm
const}$, all the components of $t''$ have $Lie_E t''\neq {\rm const}$, we
have
$
F(t) =F_0(t)+\sum_{k,\, l\geq 0}A_{k,l}(t'')^l e^{k\, t'}
$
and the series converges in some neiborghood of $t''=0$, $t'=-\infty$. 

\noindent{\bf 2.2. K.Saito theory of primitive forms and Frobenius
structures on
universal unfoldings of quasihomogeneous singularities.} Let $f_s(x)$,
$s=(s_1, \dots, s_n)$ be the universal unfolding of a quasihomogeneous
isolated singularity $f(x)$, $x\in {\bf C}^N$, $f(0)=f'(0)=0$. Here $n$ is
the Milnor number of the singularity. The Frobenius structure on the base
$M\ni s$ of the universal unfolding can be easily constructed [BV]
using the
theory of primitive forms [Sai2]. For the 
example [DVV] of the $A_n$ singularity $f(x)=x^{n+1}$ the universal
unfolding reads
$f_s(x) =x^{n+1} + s_1 x^{n-1} + \dots + s_n$, $M={\bf C}^n\ni (s_1,
\dots, s_n)$. On the FM $e=\partial/\partial s_n$, $E=\sum (k+1)s_k
\partial/\partial s_k$, the metric has the form
$$
<\partial_{s_i}, \partial_{s_j}> =-(n+1) \res_{x=\infty} {
\partial f_s(x)/\partial s_i \, \partial f_s(x) / \partial s_j
\over
f'_s(x)}
\eqno(2.1)$$
the multiplication is defined by
$$
<\partial_{s_i}\cdot\partial_{s_j}, \partial_{s_k}>=
-(n+1)\res_{x=\infty} {\partial f_s(x)/\partial s_i \, \partial f_s(x) /
\partial s_j\, \partial f_s(x)/\partial s_k
\over f'_s(x)}.
\eqno(2.2)$$
This is a polynomial FM.
The deformed
flat coordinates are given by oscillatory integrals
$$
\tilde t_c ={1\over \sqrt{z}} \int_c e^{z\, f_s(x)} d x
\eqno(2.3)$$
Here $c$ is any
1-cycle in ${\bf C}$ that goes to infinity along the
direction ${\rm Re}\, z\, f_s(x) \to -\infty$.

\noindent{\bf 2.3. Quantum cohomology} of a $2d$-dimensional smooth
projective
variety $X$ is a Frobenius structure of the charge $d$ on a domain
$M\subset H^*(X,{\bf C})/2\pi i H^2(X,{\bf Z})$ (we assume that
$H^{odd}(X)=0$ to avoid working with supermanifolds, see [KM]). It is an
analytic perturbation in the sense of n.2.1 of the cubic for $A_0=H^*(X)$
defined by a generating function of the genus zero Gromov - Witten (GW)
invariants of $X$ [W1, W2, MS, RT,
KM, Beh].  They are defined as intersection numbers of  certain cycles on
the moduli spaces of stable maps [KM]
$$
X_{[\beta], l} := \left\{ 
  \beta :(S^2, p_1, \dots,
p_l)\to X, ~~{\rm given ~homotopy ~class}~ [\beta] \in H_2(X;{\bf
Z})\right\} .
$$
The holomorphic maps $\beta$ of the Riemann sphere $S^2$ with $l\geq 1$
distinct marked points are considered up to a holomorphic change of
parameter. The markings define evaluation maps
$
p_i: X_{[\beta], l} \to X, ~~(\beta, p_1, \dots, p_l)\mapsto \beta(p_i).
$
$$
F(t) = F_0(t) +\sum_{[\beta]\neq 0}\sum_l \big < e^{t''}\big >_{[\beta],l}
\exp
\int_{S^2}\beta^*(t'), ~~\big < e^{t}\big >_{[\beta],l}:={1\over l!} 
\int_{X_{[\beta],l}}
p_1^*(t)\wedge \dots \wedge p_l^*(t)
\eqno(2.4)$$
for $t=(t', t'')\in H^*(X)$, $t'\in H^2(X)/2\pi i H^2(X,{\bf Z})$, $t''
\in H^{*\neq 2}(X)$. This potential together with the Poincar\'e pairing
on $TM=H^*(X)$, the unity vector field $e=1\in H^0(X)$, the Euler vector
field
$
E(t) =\sum (1-q_\alpha) t^\alpha e_\alpha + c_1(X), ~~t=t^\alpha
e_\alpha, ~e_\alpha \in H^{2 q_\alpha}(X)
$
gives the needed Frobenius structure. The deformed flat coordinates are
generating functions of certain ``gravitational descendents'' [Du5], see
also 
[DW, Ho, Gi1]
$
\tilde t_\alpha(t;z) =
\sum_{p=0}^\infty \sum_{[\beta],l} 
\big < 
z^{\mu +p} 
z^{c_1(X)} 
\tau_p(e_\alpha)
\otimes 1 
\otimes e^{t''}
\big >_{[\beta],l}e^{\int_{S^2} \beta^*(t')}
$,
$\alpha=1, \dots, n={\rm dim}\, H^*(X)$, $\mu(e_\alpha) =(q_\alpha-d/2)
e_\alpha$,
The definition of the descendents $\big < \tau_{p_1}(a_1)\otimes
\tau_{p_2}(a_2)\otimes \dots \otimes \tau_{p_l}(a_l)\big >_{[\beta],l}$
see in  [W2], [KM].
The definition of GW invariants can be extended on a certain class of
compact symplectic varieties $X$ using Gromov's theory [Gr] of 
pseudoholomorphic curves, see [W2, MS, RT].
\smallskip
\noindent{\bf 3. Classification of semisimple FMs.} 

\noindent{\bf 3.1. Definition.} A point $t\in M$ is called {\it
semisimple} if the
algebra on $T_tM$ is semisimple. A connected FM $M$ is
called
semisimple if it has at least one semisimple point.
Classification of semisimple FMs can be reduced, by a
nonlinear change of coordinates, to a system of ordinary differential
equations. First we will describe these new coordinates.
 
\noindent{\bf 3.2.
Canonical coordinates} on a semisimple FM. Denote $u_1(t)$, \dots, $u_n(t)$ 
the roots of the characteristic
polynomial of the operator of multiplication by the Euler vector field
$E(t)$ ($n={\rm dim}\, M$). Denote $M^0\subset M$ the open subset where
all the roots are pairwise distinct.
It turns out [Du2] that
the functions $u_1(t)$, \dots, $u_n(t)$
are
independent local coordinates on $M^0\neq \emptyset$. In these coordinates
$
\partial_i\cdot \partial_j =\delta_{ij} \partial_i, ~~{\rm
where}~\partial_i:=\partial/\partial u_i,
$
and $E=\sum_i u_i \partial_i$.
The local coordinates $u_1$, \dots, $u_n$ on $M^0$ are called {\it
canonical}. 

\noindent{\bf 3.3. Deformed flat connection in the canonical coordinates
and
isomonodromy deformations.} Staying in a small ball on $M^0$, let us order
the canonical coordinates and choose the signs of the square roots
$
\psi_{i1}:=\sqrt{<\partial_i, \partial_i>}, ~~i=1, \dots, n.
$
The orthonormal frame of the normalized idempotents $\partial_i$
establishes a local trivialization of the tangent bundle $TM^0$. The
deformed flat connection $\tilde \nabla$ in $TM^0$ is recasted into the
following flat connection in the trivial bundle $M^0\times {\bf C}\times
{\bf C}^n$
$$
\tilde\nabla_i = \partial_i -z\, E_i -V_i, ~~~~
\tilde\nabla_{d/dz} = \partial_z -U - z^{-1}V,
\eqno(3.1)$$
other components are obvious. Here the  $n\times n$ matrices $E_i$, $U$,
$V=(V_{ij})$ read
$
(E_i)_{kl}=\delta_{ik}\delta_{il}, ~~U={\rm diag}\, (u_1, \dots, u_n),
~V=\Psi \, \mu\, \Psi^{-1}=-V^T
$
where the matrix $\Psi=(\psi_{i\alpha})$ satisfying $\Psi^T\Psi=\eta$ is
defined by
$
\psi_{i\alpha} := 
\psi_{i1}^{-1}
\partial t_\alpha/ \partial u_i,
~~i, \, \alpha=1, \dots, n
$.
The skew-symmetric matrices $V_i$ are determined by the equations
$[U,V_i]=[E_i,V]$.

Flatness of the connection (3.1) reads as the system of commuting
time-dependent Hamiltonian flows on the Lie algebra $so(n)\ni V$ equipped
with the standard linear Poisson bracket
$$
\partial_iV=\left\{ V,H_i(V;u)\right\} ,~~i=1, \dots, n
\eqno(3.2)$$
with the quadratic Hamiltonians
$
H_i(V;u) ={1\over 2} \sum_{j\neq i} {V_{ij}^2\over u_i-u_j}, ~i=1, \dots,
n.
$
For the first nontrivial case $n=3$ (3.2) can be reduced to a particular
case
of the classical Painlev\'e-VI equation.
The monodromy of the operator $\tilde\nabla_{d/dz}$ (i.e., the monodromy
at the origin, the Stokes matrix, and the central connection matrix, see
definitions in [Du3, Du5]) does not change with small variations of a
point $u=
(u_1, \dots, u_n)\in M$. 

\noindent{\bf 3.4. Parametrization of semisimple FMs by
monodromy
data of the deformed flat connection.} 
We now reduce the above system of nonlinear differential equations to a
linear boundary value problem of the theory of analytic functions. First
we will describe the set of parameters of the boundary value problem.

\noindent{\bf 3.4.1.Monodromy at the origin} (defined also for 
nonsemisimple FM)
consists of

- a linear $n$-dimensional space ${\cal V}$ with a symmetric nondegenerate
bilinear form $<~,~>$, a skew-symmetric linear operator $\mu:{\cal V} \to
{\cal V}$, $<\mu(a), b>+ <a,\mu(b)>=0$, and a marked eigenvector $e_1$ of
$\mu$, $\mu(e_1) = -d/2 \, e_1$. In main examples the operator $\mu$ will
be diagonalizable.

-  A linear operator $R: {\cal V}\to {\cal V}$ satisfying the following
properties: (1) $R=R_1+R_2+\dots$ where $R_k( {\cal V}_\lambda )\subset
{\cal V}_{\lambda+k}$ for the root decomposition of ${\cal V}
=\oplus_\lambda {\cal V}_\lambda$, $\mu(v_\lambda)=\lambda v_\lambda$
for $v_\lambda \in {\cal V}_\lambda$. (2) $\left\{ Rx, y\right\}
+\left\{ x, Ry\right\}=0$ for any $x,\, y\in {\cal V}$ where
$
\left\{ x,y\right\} := \left< e^{\pi i \mu} x, y\right>.
$


\noindent{\bf 3.4.2. Stokes matrix} is an arbitrary $n\times n$ upper
triangular
matrix $S=(s_{ij})$ with $s_{ii}=1$, $i=1, \dots, n$. We treat it as a
bilinear form $<a,b>_S:=a^T S b$, $a, b\in {\bf C}^n$.

\noindent{\bf 3.4.3. Central connection matrix} is an isomorphism $C: {\bf
C}^n
\to {\cal V}$ satisfying
$
<a,b>_S =<Ca, e^{\pi i \mu} e^{\pi i R} Cb> ~{\rm for~any}~a, \, b\in {\bf
C}^n
$.
The matrices $S$ and $C$ are defined up to a transformation
$S\mapsto DSD, ~~C\mapsto CD, ~~D={\rm diag}\, (\pm 1, \dots, \pm 1).
$

\noindent{\bf 3.4.4. Riemann - Hilbert boundary value problem} (RH
b.v.p.). Let us
fix a radius $R>0$ and an argument $0\le \varphi < 2\pi$. Denote $\ell
=\ell_+ \cup \ell_-$ the oriented line $\ell_+ =\{ z \, | \, {\rm arg}\,
z=\varphi\}$, $\ell_-=\{ z\, | \, {\rm arg} \, z=\varphi +\pi\}$. It 
divides the complex $z$-plane into two halfplanes $\Pi_\rr$ and $\Pi_\ll$.
For a given $u=(u_1, \dots, u_n)$ with $u_i\neq u_j$ for $i\neq j$ and for
given monodromy data we are looking for: (1) $n\times n$ matrix-valued
functions $\Phi_\rr(z)$, $\Phi_\ll(z)$ analytic for $|z|>R$ and $z\in
\Pi_\rr$ and
$z\in \Pi_\ll$ resp., continuous up to the boundaries $|z|=R$ or $z\in
\ell$ and satisfying
$
\Phi_{\rr / \ll}(z) =1 +O\left(1/ z\right) ~~{\rm for }~ |z|\to
\infty
$
within the correspondent half-plane $\Pi_{\rr / \ll}$;
(2) $n\times n$ matrix-valued function $\Phi_0(z)$ (with values in
$Hom( {\cal V}, {\bf C}^n)$) analytic for $|z|<R$ and continuous up to the
boundary $|z|=R$, such that $\det \Phi_0(0)\neq 0$.
The boundary values of the functions must satisfy
$$
\Phi_\rr (z) e^{zU} = \Phi_\ll (z) e^{zU} S ~~ {\rm for} ~z\in \ell_+,
~|z|>R;
$$
$$
\Phi_\rr (z) e^{zU} = \Phi_\ll (z) e^{zU} S^T ~~ {\rm for}~ z\in \ell_-,
~|z|>R;
$$
$$
\Phi_\rr(z) e^{zU} =\Phi_0(z) z^\mu z^R C ~~ {\rm for }~|z|=R, ~z\in
\Pi_\rr;
$$
$$
\Phi_\ll (z) e^{zU}S = \Phi_0(z) z^\mu z^R C ~~ {\rm for} ~|z|=R, ~z\in
\Pi_\ll.
$$
The branchcut in the definition of the multivalued functions $z^\mu$ and
$z^R$ is chosen along $\ell_-$.
For solvability of the above RH b.v.p. we have also to require the complex
numbers $u_1$, \dots, $u_n$ to be ordered in such a way, depending on
$\varphi$, that 
$$
{\cal R}_{jk}:=\{ z=-ir(\bar u_j-\bar u_k) | r\geq 0\} \subset \Pi_\ll
~~{\rm for }~ j<k.
\eqno(3.3)$$
Denote ${\cal U}(\varphi)\subset {\bf C}^n$ the set of all points $u=(u_1,
\dots, u_n)$ with $u_i\neq u_j$ for $i\neq j$ satisfying (3.3). Let ${\cal
U}_0(\varphi)$ be the subset of points $u\in {\cal U}(\varphi)$ such that:
(1) the RH b.v.p. is solvable and (2) all the coordinates of the vector
$\Phi_0(0)e_1$ are distinct from zero. It can be shown (cf. [Mi],
[Mal]) that the solution $\Phi_{\rr /\ll}=\Phi_{\rr /\ll}(z;u)$,
$\Phi_0=\Phi_0(z;u)$ of the RH b.v.p depends analytically on $u\in {\cal
U}_0(\varphi)$. Let
$
\Phi_0(z;u)=\sum_{p=0}^\infty \phi_p(u) z^p.
$
Denote (only here) $(~,~)$ the standard sum of squares quadratic form on
${\bf C}^n$. Choose a basis $e_1$, $e_2$, \dots, $e_n$ of eigenvectors of
$\mu$, $\mu(e_\alpha)=\mu_\alpha e_\alpha$, $\mu_1=-d/2$, and put
$\eta_{\alpha\beta}:=<e_\alpha, e_\beta>$,
$(\eta^{\alpha\beta}):=(\eta_{\alpha\beta})^{-1}$.
\smallskip
{\bf Theorem 1} [Du2, Du3, Du5]. {\it The formulae
$$
t_\alpha(u) =(\phi_0(u) e_\alpha, \phi_1(u) e_1),  ~~
t^\alpha=\eta^{\alpha\beta}t_\beta,  ~\alpha=1, \dots, n,
$$
$$
F=1/ 2 \left[ (\phi_0t, \phi_1t)
-2(\phi_0 t, \phi_1e_1)
+(\phi_1e_1, \phi_2e_1) -(\phi_3e_1,\phi_0e_1)\right]
$$
$$
E(t)= \sum_{\alpha=1}^n (1+\mu_1-\mu_\alpha)t^\alpha \partial_\alpha
+\sum_\alpha (R_1)_1^\alpha \partial_\alpha
$$
define on ${\cal U}_0(\varphi)$ a structure of a semisimple FM
$Fr({\cal V}, <~,~>, \mu, e_1, R, S, C)$. 
Any semisimple FM locally has such a form.}

\noindent{\bf 3.5. Remark.} The columns of the matrices $\Phi_0(z;u)z^\mu
z^R$ and $\Phi_\rr(z;u)e^{zU}$ correspond to two different bases in the
space of deformed flat coordinates. The first basis is a deformation,
$z\to 0$,
of the original flat coordinates, $\tilde t =[t+O(z)] z^\mu z^R$. The
second one, defined only in the semisimple case, corresponds to a system
of deformed flat coordinates given by oscillatory integrals (see (2.3) and
Section 6 below).

\noindent{\bf 3.6. Global structure of semisimple FMs and
action
of the braid group on the monodromy data.} Let $B_n$ be the group of
braids with $n$ strands. We will glue globally the FM from
the charts described in n.3.4 with different $S$ and $C$. So, for brevity,
we redenote here  the charts $Fr({\cal V}, <~,~>, \mu, e_1, R, S, C)=:
Fr(S,C)$. The charts will be labelled by braids $\sigma\in B_n$. By
definition in the chart $Fr(S^\sigma, C^\sigma)$ the functions
$t^\alpha(u)$, $F(u)$ are obtained as the result of analytic continuation
from $Fr(S,C)$
along the braid $\sigma$. The action $S\mapsto S^\sigma$, $C\mapsto
C^\sigma$ of the standard generators $\sigma_1$, \dots, $\sigma_{n-1}$
of $B_n$  is given by
$
S^{\sigma_i} = KSK,  ~~C^{\sigma_i} =CK
$
where the only nonzero entries of the matrix $K= K^{(i)}(S)$ are
$
K_{kk}=1, ~k=1, \dots, n, ~k\neq i, \, i+1, ~~K_{i,i+1}=K_{i+1,i}=1, ~~
K_{i,i}=-s_{i,i+1}
$.
Let $B_n(S,C)\subset B_n$ be the subgroup of all braids $\sigma$ such that
$
S^\sigma =DSD, ~~C^\sigma = CD, ~~D={\rm diag} (\pm 1, \dots, \pm 1).
$

{\bf Theorem 2} [Du3, Du5]. {\it Any semisimple FM has the form
$
M=\cup_{\sigma\in B_n/B_n(S,C)}
 Fr({\cal V}, <~,~>, \mu, e_1, R, S^\sigma, C^\sigma)
$
where the glueing of the charts is given by the above action of $B_n$.}

\noindent{\bf 3.7. Tau-function of the isomonodromy deformation and
elliptic GW
invariants.} Like in n.2.2, the genus $g$ GW invariants can be defined in
terms of the intersection theory on the moduli space $X_{[\beta], l}(g)$ 
of stable
maps $\beta: C_g\to X$ of curves of genus $g$ with markings [KM, Beh]. It
turns
out that, assuming semisimplicity of quantum cohomology of $X$, the
elliptic (i.e., of $g=1$) GW invariants can still be expressed via
isomonodromy deformations. To this end we define, following [JM], the
$\tau$-function $\tau(u_1, \dots, u_n)$ of a solution $V(u)$ of the system
(3.2) by the quadrature of a closed 1-form
$
d\log \tau = \sum_{i=1}^n H_i(V(u);u) du_i
$.
We also define $G$-function of the FM by
$
G=\log ({\tau/ J^{1/24}})
$
where 
$
J=\det \left( \partial t^\alpha /\partial u_i\right) =\pm \prod_{i=1}^n
\psi_{i1}(u).
$

{\bf Theorem 3} [DZ2]. {\it For an arbitrary semisimple FM the
$G$-function is the unique, up to an
additive constant, solution to the system of 
[Ge] for the generating function of elliptic GW invariants
satisfying 
$
Lie_e G=0, ~~Lie_E G={\rm const}.
$
}

\noindent{\bf 3.8. Problem of selection of semisimple FMs
with good
analytic properties} of n.2.1 is still open. Experiments for small $n$
[Du3]
show
that such solutions are rare exceptions among all semisimple FMs.
Analyticity of the $G$-function near the point $t'=-\infty$,
$t''=0$ imposes further restrictions on $M$ [DZ2]. 
To solve the problem one is to study the behaviour of solutions of the RH
b.v.p. in the limits when two or more among the canonical coordinates
merge. At the point  $t'=-\infty$,
$t''=0$ all $u_1=\dots =u_n=0$.

\smallskip
\noindent{\bf 4. Examples of monodromy data.} 

\noindent{\bf 4.1. Universal unfoldings of isolated singularities.} The
subspace
$M_0\subset M$ consists of the parameters $s$ for which the versal
deformation $f_s(x)$ has $n={\rm dim}\,M$ distinct critical values
$u_1(s)$, \dots, $u_n(s)$. These will be our canonical coordinates. The
monodromy at the origin is the classical monodromy operator [AGV] of the
singularity, the Stokes matrix coincides with the matrix of the variation
operator in the Gabrielov's distinguished basis of vanishing cycles (see
[AGV]; we may assume that ${\rm dim}\, x\equiv 1 \, ({\rm mod} 4)$).

\noindent{\bf 4.2. Quantum cohomology of Fano varieties.}
The following two questions are to be answered in order to apply the above
technique to the quantum cohomology of a variety $X$.

{\bf Problem 1.} When do the generating series (2.4) converge?

{\bf Problem 2.} For which $X$ the quantum
cohomology of $X$ is semisimple? 

Hopefully, in the semisimple case the convergence can be proved on the
basis of the differential equations of n.3. To our opinion the problem 2
is more deep. A necessary condition to have a semisimple quantum
cohomology
is that $X$ must be a Fano variety. It was conjectured to be also a
sufficient condition [TX], [Man1]. We analyze below one example and
suggest some more modest conjecture describing also a part of the
monodromy data.

\noindent{\bf 4.2.1. Quantum cohomology of projective spaces.} For $X={\bf
P}^d$:
(1) the monodromy at the origin is given by the bilinear form $<e_\alpha,
e_\beta>=\delta_{\alpha+\beta, d+2}$
in
$H^*(X)={\cal V} ={\rm span}\, (e_1, \dots, e_{d+1})$, the matrix 
$
\mu=1/2 {\rm diag}\, \left( -d, 1-d, \dots, d-1,
d\right)
$
and $R$ is the matrix of multiplication by the first Chern class
$
R=R_1 =c_1(X), ~~R\, e_\alpha =(d+1) e_{\alpha+1} ~{\rm for}~ \alpha\leq
d, ~R\,
e_{d+1}=0.
$
With obvious modifications these formulae work also for any variety $X$
with $H^{\rm odd}(X)=0$ (see [Du3]).
(2) The Stokes matrix $S=(s_{ij})$ has the form
$$
s_{ij}=\left(\matrix{ d+1\cr j-i\cr}\right) ~{\rm for}~ i\leq j,
~~s_{ij}=0 ~{\rm for }~ i>j.
\eqno(4.1)$$
This form of Stokes matrix was conjectured in [CV], [Zas] but, to our
knowledge, it was proved only in [Du5] for $d=2$ and in [Guz] for any $d$.
(3) The central connection matrix $C$ has the form $C=C' C''$,
$C'=\left( {C'}_\beta^\alpha\right)$, $C''=\left({C''}_j^\beta\right)$ 
where
$
{C''}_j^\beta =[2 \pi i (j-1)]^{\beta -1}/ (\alpha-1)!, ~~j,\, \beta =1,
\dots, d+1$,
${C'}_\beta^\alpha ={(-1)^{d+1}\over (2\pi)^{d+1\over 2} i^{\bar d}}
\cases{A_{\alpha-\beta} (d), ~\alpha \geq \beta\cr
0, ~\alpha < \beta\cr}
$
with $\bar d=1$ for $d=$even and $\bar d=0$ for $d=$odd
where the numbers $A_0(d)=1$, $A_1(d)$, \dots, $A_d(d)$ are defined
from the Laurent expansion for $x\to 0$:
$
1/ x^{d+1} +A_1(d)/ x^d + \dots + A_d(d)/ x + O(1) =
(-1)^{d+1}\Gamma^{d+1} (-x) e^{-\pi i \bar d x} 
$.
Observe that (4.1) is the Gram matrix of the bilinear form $\chi(E,F)
:= \sum_k (-1)^k {\rm dim}\, Ext^k(E,F)$ in the basis 
given by a
particular full system $E_j ={\cal O}(j-1)$, $j=1, \dots, d+1$ of
exceptional objects in the derived category $Der^b(Coh({\bf P}^d))$ of
coherent sheaves on ${\bf P}^d$ [Rud]. The columns of the matrix
$C''$ are the components of the Chern character ${\rm ch}\, (E_j) =e^{2
\pi i c_1 (E_j)}$, $j=1, \dots, d+1$. The geometrical meaning of the
matrix $C'$ remains unclear. 
In other charts of the FM $S^\sigma$ and $C^\sigma =C'
{C''}^\sigma$, $\sigma\in B_n$, have the same structure for another full
system $E_1^\sigma$, \dots, $E_{d+1}^\sigma \in Der^b(Coh({\bf P}^d))$
of exceptional objects, where the action of the braid group
$(E_1, \dots, E_{d+1}) \mapsto (E_1^\sigma, \dots, E_{d+1}^\sigma)$ is
described in [Rud]. Warning: the points of the FM
corresponding to the restricted quantum cohomology [MM], where $t\in
H^2({\bf
P}^d)$, do not belong to the chart $Fr(S,C)$ with the matrices $S$ and $C$
as above!

\noindent{\bf 4.2.2. Conjecture.} We say that a Fano variety  $X$ is {\it
good} if
$Der^b(Coh(X))$ admits, in the sense of [BP], a full system of exceptional
objects $E_1$, \dots, 
$E_n$, $n={\rm dim}\, H^*(X)$. Our conjecture is that (1) the quantum
cohomology of $X$ is semisimple {\it iff} $X$ is a good Fano variety; (2)
the Stokes matrix $S=(s_{ij})$ is equal to $s_{ij}=\chi(E_i,E_j)$, $i, \,
j=1, \dots, n$; (3) the central connection matrix has the form $C=C' C''$
when the columns of $C''$ are the components of ${\rm ch}\, (E_j) \in
H^*(X)$ and $C':H^*(X)\to H^*(X)$ is some operator satisfying $C'
(c_1(X)a) =c_1(X) C'(a)$ for any $a\in H^*(X)$.

For $X={\bf P}^d$ the validity of the conjecture follows from n.4.2.1
above. The conjecture probably can be derived from more general conjecture
[Kon] about equivalence of $Der^b(Coh(X))$ to the Fukaya category of the
mirror pair $X^*$ of $X$. According to it (see also [EHX, Gi1, Gi2])
the
basis of horizontal sections of $\tilde\nabla$ corresponding to the
columns of $\Phi_\rr (z; u) e^{zU}$ coincides with the oscillatory 
integrals of the Fukaya category of $X^*$. However, we do not know who
is the first factor $C'$ of the connection matrix in this general setting.  
\smallskip
\noindent{\bf 5. Intersection form} of a FM is a bilinear
symmetric
pairing on $T^*M$ defined by $(\omega_1, \omega_2)|_t := i_{E(t)}
(\omega_1\cdot \omega_2)$, $\omega_1, \, \omega_2 \in T_t^*M$. {\it
Discriminant} is the locus $\Sigma=\{ t\in M \, | \, \det (~,~)_t =0 \}$.
On $M\setminus \Sigma$ the inverse to $(~,~)$ determines a flat metric
and, thus, a local isometry $\pi : \left( M\setminus \Sigma,
(~,~)^{-1}\right) \to {\bf C}^n$ where ${\bf C}^n$ is equipped with a
constant complex Euclidean metric $(~,~)_0$. This local isometry is called
{\it period mapping} (our terminology copies that of the singularity
theory
where the geometrical structures with the same names live on the bases of
universal unfoldings, see [AGV]). The image $\pi (\Sigma)$ is a collection
of nonisotropic hyperplanes in ${\bf C}^n$. Multivaluedness of $\pi$ is
described by the {\it monodromy representation} $\pi_1 (M\setminus \Sigma)
\to Iso\left({\bf C}^n, (~,~)_0 \right)$ (for $d\neq 1$ to the orthogonal
group $O\left( {\bf C}^n, (~,~)_0\right)$). The image $W(M)$ of the
representation is called {\it monodromy group} of the FM
$M$. In the semisimple case it is always an extension of a reflection
group (see details in [Du5]). Our hope is that, for a semisimple
FM $M$ with good analytic properties, the monodromy group
acts discretely in some domain $\Omega\subset {\bf C}^n$, and $M$ is
identified with a branched covering of the quotient $\Omega / W(M)$.

\noindent{\bf 5.1. Examples} of a FM with $W(M)$ = finite
irreducible Coxeter
group $W$ acting in ${\bf R}^n$ [Du3]. 
%
%
%
%
%
%
These are polynomial FMs, $M={\bf C}^n/W$, constructed in terms of the
theory of invariant polynomials of $W$.
Conjecturally, all polynomial semisimple FMs are
equivalent to the above and to their direct sums. 

This construction was generalized in [DZ1] to certain extensions of affine
Weyl groups and in [Ber] to Jacobi groups of the types $A_n$, $B_n$,
$G_2$. For the quantum cohomology of
${\bf P}^2$ the monodromy group is isomorphic to $PSL_2({\bf Z})\times \{
\pm 1\}$ [Du5].
\smallskip
\noindent{\bf 6. Mirror construction} represents certain system of
deformed flat
coordinates on a semisimple FM by oscillatory integrals
$
I_j (u; z) ={1\over \sqrt{z}} \int_{Z_j} e^{z \lambda (p; u)} dp,
~~\tilde\nabla I_j(u;z)=0, ~~j=1, \dots, n
$
having the phase function $\lambda(p;u)$ depending on the
parameters $u=(u_1, \dots, u_n)$ defined on a certain family of open
Riemann surfaces ${\cal R}_u\ni p$ realized as a finite-sheeted branched 
covering $\lambda: {\cal R}_u \to D\subset {\bf C}$ over a domain in the
complex plain. The ramification points of ${\cal R}_u$, i.e., the critical
values of the phase function,  
are $u_1$, \dots,
$u_n$. The 1-cycles $Z_1$, \dots, $Z_n$ on ${\cal R}_u$ go to infinity in
a way that guarantees the convergence of the integrals. The function
$\lambda(p;u)$ satisfies an important property: for any two critical
points $p_i^{1, \, 2}\in {\cal R}_u$ with the same critical value $u_i$
the equality $d^2 \lambda (p_i^1; u) /dp^2 = d^2 \lambda (p_i^2; u)/dp^2$
must hold true. The metric $<~,~>$  and
the trilinear form $c(a_1, a_2, a_3) := <a_1\cdot a_2, a_3>$ are given by
the residue formulae similar to (2.1), (2.2).
The solutions $p=p(u;\lambda)$ of the
equation $\lambda(p;u)=\lambda$ are the flat coordinates of the
{\it flat pencil of the metrics} $(~,~)-\lambda<~,~>$ on $T^*M$ [Du3-Du5].

For the case when generically there is a unique critical point $p_i$ over
$u_i$ for each $i$ and ${\cal R}_u$ can be compactfied to a Riemann
surface of a finite genus $g$, we arrive at the Hurwitz spaces of branched
coverings [Du1, Du3].

The construction of the Riemann surfaces ${\cal R}_u$, of the phase
function $\lambda (p;u)$ and of the cycles $Z_1$, \dots, $Z_n$ is given in
[Du5] by universal formulae assuming $\det (S+S^T)\neq 0$. In the quantum
cohomology of a $d$-fold $X$ the last condition is valid for $d$ = even.
For $d$ = odd one has $\det (S-S^T)\neq 0$. In this case one can represent
the deformed flat coordinates  by oscillatory integrals with the phase
function $\lambda(p,q;u)=\nu (p;u) + q^2$ depending on two variables $p$,
$q$. The details will be published elsewhere.
\smallskip
\noindent{\bf 7. Gravitational descendents} is a physical name for
intersection
numbers $<\tau_{m_1}(a_1)\otimes \dots \otimes\tau_{m_l}(a_l)>$ of the
pull-back
cocycles $p_1^*(a_1)$, \dots, $p_l^*(a_l)$ with the Mumford - Morita -
Miller cocycles $\psi_1^{m_1}$, \dots, $\psi_l^{m_l}\in H^*
(X_{[\beta],l})$ [W2], [DW], [KM]. We will describe first their genus
$g=0$
generating function
$
{\cal F}_0(T) = \sum _{[\beta]} \sum_l \big<
e^{\sum_{\alpha=1}^n\sum_{p=0}^\infty T^{\alpha,p}
\tau_p(e_\alpha)}\big>_{ [\beta], l, g=0}.
$
Here $T=\left(T^{\alpha,p}\right)$ are indeterminates (the coordinates on
the ``big phase
space'', according to the physical terminology). This function has the
form 
$
{\cal F}_0(T)=1/ 2 \sum \Omega_{\alpha,p;\, \beta,q} \left(
t(T)\right) \tilde T^{\alpha,p}\tilde T^{\beta,q}
$
where
$
\tilde T^{\alpha,p} = T^{\alpha,p} ~{\rm for }~ (\alpha,p)\neq (1,1), ~
\tilde T^{1,1} = T^{1,1}-1,
$
the functions $\Omega_{\alpha,p;\, \beta,q}(t)$ on $M$ are the
coefficients of the expansion of the matrix valued function
$
\Omega_{\alpha\beta}(z,w;t) := (z+w)^{-1} \left[ \left(\Phi_0^T(w;t)
\Phi_0(z;t)\right)_{\alpha\beta} -\eta_{\alpha\beta}\right]
=\sum_{p,\, q \geq 0} \Omega_{\alpha,p;\, \beta,q} (t) z^p w^q,
$
the vector function $t(T) = \left( t^1(T), \dots, t^n(T)\right)$ 
$$
t^\alpha (T) = T^{\alpha,0} + \sum_{q>0} T^{\beta,q} \nabla^\alpha 
\Omega_{\beta,q;\, 1,0} (t)|_{t^\alpha=T^{\alpha,0}} + \dots
\eqno(7.1)$$
is
defined as the unique solution of the following fixed point equation
$
t=\nabla \sum_{\alpha,p} T^{\alpha,p}\Omega_{\alpha,p;\, 1,0}(t).
$

The generating function ${\cal F}_1(T)$ of the genus $g=1$ descendents 
has the form [DZ2], [DW], [Ge]
$
{\cal F}_1(T) =\left[ G(t) +{1\over 24} \log \det M_{\alpha\beta}(t, \dot
t)\right] _{t=t(T), \, \dot t= \partial _{T^{1,0}} t(T)}
$
where $G(t)$ is the $G$-function of the FM, the matrix
$M_{\alpha\beta}(t, \dot t)$ has the form $M_{\alpha\beta}(t, \dot t) = 
\partial_\alpha \partial_\beta \partial_\gamma F(t)\dot t^\gamma$, the
vector
function $t(T)$ is the same as above. The structure of the genus $g=2$
corrections is still unclear, although there are some interesting
conjectures [EX] related, in the case of quantum cohomology, to the
Virasoro constraints for the full partition function
$$
Z(T; \varepsilon) =\exp \sum_{g=0}^\infty \varepsilon^{2g-2} {\cal
F}_g(T)
\eqno(7.2)$$
$\varepsilon$ is a formal small parameter called {\it string coupling
constant}.
\smallskip
\noindent{\bf 8. Integrable hierarchies} of PDEs of the KdV type and
FMs.
The idea that FMs may serve as moduli of
integrable hierarchies of evolutionary equations (see [W2], [Du2], [Du3])
is based
on

(1) the theorem of Kontsevich - Witten identifying the partition function
(7.2) in the case $X$ = point as the tau-function of a particular solution
of the KdV hierarchy.

(2) The construction [Du2, Du3] of bihamiltonian integrable hierarchy of
the
Whitham type
$
\partial_{T^{\alpha,p}} t =\left\{ t(X),
H_{\alpha,p}\right\}_1
=   K^{(0)}_{\alpha,p}(t,t_X)
$
(the vector function in the r.h.s. depends linearly on the derivatives
$t_X$) such that the full genus zero partition function is the
tau-function of a particular solution (7.1) to the hierarchy. The
solution is specified by the symmetry constraint $t_X -\sum
T^{\alpha,p} \partial_{T^{\alpha,p-1}}t =1$. The
phase space of the hierarchy  is the loop space ${\cal L}(M) = \left\{
\left( t^1(X), \dots, t^n(X)\right)\, |\, X\in S^1\right\}$, the first
Hamiltonian structure is
$
\left\{ t^\alpha(X), t^\beta(Y)\right\}_1 = \eta^{\alpha\beta}
\delta'(X-Y),
$
the second one $\{ ~,~\}_2$ is determined [Du3] by the flat metric $(~,~)$
according to the
general scheme of [DN]. The Hamiltonians are $H_{\alpha,p}=\int
\Omega_{\alpha,p;\, 1,0}(t)\, dX$. Actually, any linear combination
$\{~,~\}_2 -\lambda\{~,~\}_1$ with an arbitrary $\lambda$ is again a
Poisson bracket on the loop space since $(~,~)$ and $<~,~>$ form {\it
a flat pencil of metrics} on $T^*M$ [Du3, Du4]
(this {\it bihamiltonian property}
is a manifestation of integrability of the hierarchy, see [Mag], [Du4]). 

What we want to construct is a deformation of the hierarchy of the
form
$
\partial_{T^{\alpha,p}} t = K^{(0)}_{\alpha,p}(t,t_X) + 
\sum_{g \geq 1} \varepsilon^{2g} K^{(g)}_{\alpha,p}(t,t_X, \dots,
t^{(2g+1)})
$
where $K^{(g)}_{\alpha,p}$ are some vector valued polynomials in $t_X$,
\dots, $t^{(2g+1)}$ with the coefficients depending on $t\in M$.
All the equations of the hierachy must commute pairwise. The full
partition function must be the tau-function of a particular solution
to the hierachy. The first $g=1$ correction for an arbitrary semisimple
FM was constructed in [DZ2]. Its bihamiltonian structure
is described, for $d\neq 1$, by a nonlinear deformation of the Virasoro
algebra with the central charge
$
c=6 \varepsilon^2 (1-d)^{-2} [n- 4 {\rm tr}\, \mu^2].
$
For the FMs corresponding to the $ADE$
Coxeter groups this formula gives the known result [FL] for the central
charge
of the classical
$W$-algebra of the $ADE$-type $c=12 \varepsilon^2
\rho^2$, where $\rho$ is equal to the half of the sum of positive roots
of the corresponding root system. 

More recently it has been proved [DZ3] for a semisimple FM
that the partition function (7.2) is annihilated, within the genus one
approximation, by half of a remarkable Virasoro algebra described in terms
of the monodromy data of the FM.
\smallskip
{\bf Acknowledgments.} I am grateful to D.Orlov for helpful discussion of
derived categories, and to S.Barannikov and M.Kontsevich for fruitful
conversations.
\smallskip
\centerline{\bf References}
\eightrm
\medskip
\item{[AGV]} Arnol'd, V.I., Gusein-Zade, S.M. and Varchenko, A.N.:
Singularities
of Differentiable Maps, volumes I, II, Birkh\"auser, Boston-Basel-Berlin,
1988.

\item{[Beh]} Behrend, K.: Gromov - Witten invariants in algebraic
geometry,
{\eightit  Inv. Math.} {\eightbf  124} (1997) 601 - 627.

\item{[Ber]} Bertola, M.: Jacobi groups, Hurwitz spaces, and Frobenius
structures, Preprint SISSA 69/98/FM.


\item{[BV]} Blok, B. and Varchenko, A.: Topological conformal 
field theories
and the flat coordinates, {\eightit  Int. J. Mod. Phys.} {\eightbf  A7} (1992) 1467.

\item{[BP]} Bondal, A.I. and Polishchuk, A.E.: Homological properties
of associative algebras: the method of helices, {\eightit Russ. Acad. Sci.
Izv. Math.} {\eightbf 42} (1994) 219 - 260.

\item{[CV]} Cecotti, S. and  Vafa, C.: On classification of $N=2$ 
supersymmetric
theories, {\eightit  Comm. Math. Phys.} {\eightbf  158} (1993), 569-644.

\item{[DVV]} Dijkgraaf, R., Verlinde, E. and Verlinde, H.: Topological
strings in $d<1$, {\eightit  Nucl. Phys.}
{\eightbf  B 352} (1991) 59.

\item{[DW]} Dijkgraaf,  R., and Witten, E.: Mean field theory, topological
field theory, and multimatrix models,
{\eightit  Nucl. Phys.} {\eightbf  B 342} 
(1990) 486-522.

\item{[Du1]} Dubrovin, B.: Hamiltonian formalism of Whitham-type
hierarchies
and topological Landau - Ginsburg models, {\eightit  Comm. Math. Phys.}
{\eightbf  145} (1992) 195 - 207.

\item{[Du2]} ---: Integrable systems in topological field theory,
{\eightit  Nucl. Phys.} {\eightbf  B 379} (1992) 627 - 689.
\item{[Du3]} ---: Geometry of 2D topological field theories,
In: ``Integrable Systems and Quantum Groups'', Eds. M.Francaviglia, S.Greco,
Springer Lecture Notes in Math. {\eightbf  1620} (1996) 120 - 348.

\item{[Du4]} ---: Flat pencils of metrics and Frobenius manifolds,
math.DG/9803106, to appear in Proceedings of 1997
Taniguchi Symposium ``Integrable Systems and Algebraic Geometry''.

\item{[Du5]} ---:  Painlev\'e transcendents in two-dimensional
topological field theory, math.AG/9803107.

\item{[DN]} ---, Novikov, S.P.: The Hamiltonian formalism of
one-dimensional systems of hydrodynamic type and the Bogoliubov - Whitham
averaging method, {\eightit  Sov. Math. Dokl.} {\eightbf  27} (1983) 665 - 669.

\item{[DZ1]} ---, Zhang, Y.: Extended affine Weyl groups and 
Frobenius manifolds, {\eightit  Compositio Math.} {\eightbf  111} (1998) 167-219.

\item{[DZ2]} ---, ---: Bihamiltonian hierarchies in 2D 
topological field theory at one-loop approximation, Preprint 
SISSA 152/97/FM, hep-th/9712232, to appear in {\eightit  Comm. Math. Phys.}.

\item{[DZ3]} ---, ---: Frobenius manifolds and Virasoro
constraints, to appear

\item{[EHX]} Eguchi, T., Hori, K., Xiong, C.-S.: Gravitational quantum
cohomology,  {\eightit  Int.J.Mod.Phys.} {\eightbf  A12} (1997) 1743-1782.

\item{[EX]} Eguchi, T., Xiong, C.-S.:  Quantum Cohomology at Higher Genus:
Topological Recursion
Relations and Virasoro Conditions,
 hep-th/9801010.

\item{[FL]} Fateev, V., Lukyanov, S.: Additional symmetries and exactly
solvable models in two-dimensional conformal field theories, Parts I, II
and III, {\eightit  Sov. Sci. Rev.} {\eightbf  A15} (1990) 1. 

\item{[Ge]} Getzler, E.:  Intersection theory on $\bar M_{1,4}$ and
elliptic
Gromov-Witten invariants, alg-geom/9612004.

\item{[Gi1]} Givental, A.B.: Stationary phase integrals, quantum Toda 
lattice, flag manifolds, and the mirror conjecture, alg-geom/9612001.

\item{[Gi2]} ---: Elliptic Gromov-Witten invariants
and generalized mirror conjecture,
math.AG/9803053.

\item{[Gr]} Gromov, M.: Pseudo-holomorphic curves in symplectic
manifolds, {\eightit  Invent. Math.} {\eightbf  82} (1985), 307.

\item{[Guz]} Guzzetti, D.: Stokes matrices and monodromy groups of
the quantum cohomology of projective space, to appear.

\item{[Ho]} Hori, K.: Constraints for topological strings in $D\geq 1$,
{\eightit  Nucl. Phys.} {\eightbf  B 439} (1995) 395 - 420.

\item{[JM]} Jimbo, M. and Miwa, T.: Monodromy preserving deformations
of linear ordinary differential equations with rational coefficients. II.
{\eightit  Physica} {\eightbf  2D} (1981) 407 - 448.


\item{[Kon]} Kontsevich. M.: Talk at Scuola Normale Superiore, Pisa, April
'98.

\item{[KM]} ---, Manin, Yu.I.: Gromov - Witten classes, quantum
cohomology and enumerative geometry, {\eightit  Comm. Math. Phys.} {\eightbf  164}
(1994) 525 - 562.

\item{[MS].} McDuff, D. and  Salamon, D.: J-holomorphic curves and quantum
cohomology,
Providence, RI., American Mathematical Society, 1994.

\item{[Mag]} Magri, F.: A simple model of the integrable Hamiltonian
systems, {\eightit  J. Math. Phys.} {\eightbf  19} (1978) 1156 - 1162.

\item{[Mal]} Malgrange, B.: \'Equations Diff\'erentielles \`a
Coefficients Polynomiaux, Birkh\"auser, 1991.

\item{[Man1]} Manin, Yu.I.: Frobenius manifolds, quantum cohomology, and
moduli spaces, Preprint MPI 96-113.

\item{[Man2]} ---: Three constructions of Frobenius manifolds:
a comparative study, math.AG/9801006.

\item{[MM]}  ---, Merkulov, S.A.: Semisimple Frobenius
(super)manifolds
and quantum cohomology of $P^r$, alg-geom/9702014.
 
\item{[Mi]} Miwa, T.: 
Painlev\'e property of monodromy preserving
equations and the analyticity of $\tau$-functions, {\eightit  Publ. RIMS}
{\eightbf  17} (1981), 703-721.

\item{[Rud]} Rudakov, A.: Integer valued bilinear forms and vector
bundles,
{\eightit  Math. USSR Sbornik} {\eightbf  66} (1989), 187 - 194. 

\item{[RT]} Ruan, Y. and Tian, G.: A mathematical theory of quantum
cohomology,
{\eightit  Math. Res. Lett.} {\eightbf  1} (1994), 269-278.

\item{[Sai1]} Saito, K.: On a linear structure of a quotient variety 
by a finite
reflection group, Preprint RIMS-288 (1979), 
{\eightit  Publ.~RIMS, Kyoto Univ.}, {\eightbf  29} (1993) 535--579.

\item{[Sai2]} ---: Period mapping associated to a primitive form,
{\eightit  Publ. RIMS} {\eightbf  19} (1983) 1231 - 1264.

\item{[SYS]} ---, Yano, T. and Sekeguchi, J.: On a certain generator
system  
of the ring of invariants of a finite reflection group,
{\eightit  Comm. in Algebra} {\eightbf  8(4)} (1980) 373 - 408.

\item{[TX]} Tian, G. and Xu, G.: On the semisimplicity of the quantum 
cohomology algebra of complete intersections, alg-geom/9611035.

\item{[W1]} Witten, E.: On the structure of the topological phase of
two-dimensional
gravity,
{\eightit  Nucl. Phys.} {\eightbf  B 340} (1990) 281-332.

\item{[W2]} ---: Two-dimensional gravity and intersection theory on
moduli
space,  
{\eightit  Surv. Diff. Geom.} {\eightbf  1}  (1991) 243-210.

\item{[Zas]} Zaslow, E.:  Solitons and Helices: The Search for a
Math-Physics Bridge, {\eightit  Comm. Math. Phys.} {\eightbf  175} (1996) 337-376.
\Addresses
\end